\documentclass[11pt,reqno]{amsart}
\newtheorem{prop}{Proposition}[section]
\newtheorem{thm}{Theorem}[section]

\newtheorem{rmk}{Remark}[section]

\newcommand{\mysection}[1]{\section{#1}\setcounter{equation}{0}}

\newfont{\bb}{msbm10 at 12pt}

\def\pf{{\textit {Proof :} }}

\def\R{\hbox{\bb R}}

\def\E{\mathcal E}
\def\P{\mathcal P}

\newcommand{\bal}{\begin{align}}      \newcommand{\eal}{\end{align}}
\newcommand{\ba}{\begin{array}}      \newcommand{\ea}{\end{array}}
\newcommand{\bc}{\begin{center}}     \newcommand{\ec}{\end{center}}
\newcommand{\be}{\begin{enumerate}}  \newcommand{\ee}{\end{enumerate}}
\newcommand{\beq}{\begin{eqnarray}}  \newcommand{\eeq}{\end{eqnarray}}
\newcommand{\beQ}{\begin{eqnarray*}} \newcommand{\eeQ}{\end{eqnarray*}}
\newcommand{\bi}{\begin{itemize}}    \newcommand{\ei}{\end{itemize}}
\newcommand{\bt}{\begin{tabular}}    \newcommand{\et}{\end{tabular}}
\newcommand{\bdm}{\begin{displaymath}} \newcommand{\edm}{\end{displaymath}}

\newcommand{\rw}{\rightarrow}

\let\pa=\partial

\def\qed{\hfill{Q.E.D.}\smallskip}

\newcommand{\ls}{\setlength{\baselineskip}{12pt}
                 \setlength{\parskip}{3mm}}

\begin{document}

\title[Proceedings of ICCM 2004, December 17-22, Hong Kong]{The Positive Mass Theorem near null infinity}

\author{Xiao Zhang}
\address{Institute of Mathematics\\Academy of Mathematics and
System Sciences\\Chinese Academy of Sciences\\Beijing 100080,
China }
\email{xzhang@amss.ac.cn}

\begin{abstract}
In this short paper, we review recent progress on the positive
mass theorem for spacelike hypersurfaces which approach to null
infinity in asymptotically flat spacetimes. We use it to prove, if
the functions $c(u, \theta, \psi)$, $d(u, \theta, \psi)$ vanish at
certain retarded time in vacuum Bondi's radiating spacetimes, then
the Bondi mass is nonnegative up to that time.
\end{abstract}


\thanks{Research is partially supported by National Natural
Science Foundation of China under grant No. 10231050 and the
Innovation Project of Chinese Academy of Sciences.}


\maketitle \pagenumbering{arabic}

\mysection{Introduction}
\ls

The definition of the total energy-momentum at spatial infinity
was given by Arnowitt-Desser-Misnerfor in asymptotically flat
spacetimes \cite{ADM}. The positivity of the ADM mass was proved
by Schoen and Yau in a nontrivial isolated physical system which
satisfies the dominant energy condition \cite{SY1, SY2, SY3}.
Later it was proved by Witten using spinors \cite{Wi}. The
positive mass theorem plays a fundamental role in general
relativity. It indicates the existence of the ground state in
gravity.

Related to the gravitational radiation, there is a positive mass
conjecture at null infinity. Physically, it is believed when
gravitational radiation occurs, the energy of the system will be
carried away by gravitational waves. Most physical systems have a
conserved stress-energy tensor with a positive timelike component.
Therefore they cannot radiate away more energy than they have
initially. However, the gravitational field does not have a
well-defined stress-energy tensor. An isolated gravitational
system with finite ADM total energy-momentum might be able to
radiate arbitrarily large amounts of energy. That it cannot happen
is known as the positive mass conjecture at null infinity.

It is not available yet to set this conjecture in a mathematical
rigorous way. In the pioneering work of Bondi, van der Burg,
Metzner and Sachs on the gravitational waves in vacuum Bondi's
radiating spacetimes, the Bondi mass associated to each null cone
is defined and it is always non-increasing with respect to the
retarded time \cite{BBM, Sa, vdB}. Therefore, the Bondi mass can
be interpreted as the total mass of the isolated physical system
measured after the loss due to the gravitational radiation up to
that time. In this case, the positive mass conjecture at null
infinity is equivalent to the positivity of Bondi mass.

The proof of the positivity of the Bondi mass was outlined by
Schoen and Yau by modifying their arguments in the proof of the
positivity of the ADM mass \cite{SY4}, and by physicists by
applying Witten's spinor method, eg, see \cite{IN, HP, AsHo, LV,
ReuT, HT}. However, it needs to fill out all mathematical details.

The paper is organized as follows: In Section 2, we review
Schoen-Yau's positive mass theorem at spatial infinity. In Section
3, we review author's positive mass theorem near null infinity. In
Section 4, we review Huang-Yau-Zhang's theorem concerning the
positivity of the Bondi mass in Bondi's radiating spacetimes.
\mysection{The positive mass theorem at spatial infinity}
\ls

A spacetime is a 4-dimensional Lorentzian manifold $({\bf L}
^{3,1}, \bf{g})$ which satisfies the Einstein field equations
 \beQ
{\bf R} _{\alpha \beta}-\frac{{\bf R}}{2}\;{\bf g} _{\alpha \beta}
={\bf T} _{\alpha \beta},
 \eeQ
where ${\bf R} _{\alpha \beta}$ is the Ricci curvature of ${\bf
g}$ and ${\bf R}$ is the scalar curvature of ${\bf g}$, ${\bf T}
_{\alpha \beta}$ is the energy-momentum tensor.

There are some exact solutions of the vacuum field equations
(i.e., ${\bf T} _{\alpha \beta} =0$) in polar coordinates $(r,
\theta, \psi)$ where $0 <r <\infty$, $0\leq \theta <\pi$, $0\leq
\psi <2\pi$:
 \bi
\item Minkowski spacetime
 \beQ
 {\bf g} _{Mink}= -dt ^2 +dr ^2+r ^2 \Big(d \theta ^2 +
 \sin ^2 \theta d \phi ^2\Big);
 \eeQ
\item Schwarzschild spacetime
 \beQ
 {\bf g} _{Sch} = -\Big(1-\frac{2m}{r}\Big)dt ^2
           +\frac{dr ^2}{1-\frac{2m}{r}}
           +r ^2 \Big(d \theta ^2 +\sin ^2 \theta d \phi ^2\Big);
 \eeQ
\item Kerr spacetime
 \beQ
 {\bf g} _{Kerr} &=& -\Big(1-\frac{2mr}{\Sigma }\Big)dt ^2
             - \frac{4mar \sin ^2 \theta }{\Sigma } dt d\phi
          +\frac{\Sigma }{\triangle} dr ^2 \\
         & &+\Sigma d\theta ^2
             +\Big(r ^2 +a ^2 +\frac{2mr a^2 \sin ^2 \theta}{\Sigma ^2}
              \Big)\sin ^2 \theta d \phi ^2
 \eeQ
where $\Sigma \equiv r ^2 + a ^2 \cos ^2 \theta $, $\triangle
\equiv r ^2 -2 m r + a^2 $.
 \ei
The parameters $m$ is the total mass, $a$ is the angular momentum
per unit mass.

Let $(M, g, h)$ be a spacelike hypersurface in ${\bf L} ^{3,1}$,
$M$ is a 3-dimensional manifold, $g$ is the Riemannian metric of
$M$ and $h$ is the second fundamental form of $M$. It is usually
called an initial data set. $(M, g, h)$ is asymptotically flat if
there is a compact set $K \subset M$ such that $M \setminus K$ is
the disjoint union of a finite number of subsets $M_1,\cdots, M
_l$ - called the ``ends'' of $M$ - each diffeomorphic to $\R ^3
\setminus B _ r$, where $B _r$ is the closed ball of radius $r$
with center at the coordinate origin. In each end, $g$, $h$
satisfy
 \beQ
g_{ij}=\delta _{ij}+O\big(\frac{1}{r}\big), \partial _k g _{ij}
=O\big(\frac{1}{r ^2}\big), \partial _l \partial _k g _{ij}
=O\big(\frac{1}{r ^3}\big), h _{ij} =O\big(\frac{1}{r ^2}\big),
\partial _k h _{ij} =O\big(\frac{1}{r ^3}\big)
 \eeQ
where $\{x ^i \}$ is the Euclidean coordinates of $\R ^3$.

The total energy $E_l$ and the total linear momentum $P_{lk}$ of
end $M_l$ are defined by
 \beQ
E_l &=& \frac{1}{16\pi}\lim _{r \rightarrow \infty} \int
_{S_{r,l}}(\partial _j g_{ij}-\partial _i g_{jj})\ast dx ^i,\\
P_{lk} &=& \frac{1}{8\pi}\lim _{r \rightarrow \infty} \int
_{S_{r,l}}(h_{ki}-g _{ki} h_{jj})\ast dx ^i,
 \eeQ
where $S _{r,l}$ is the sphere of radius $r$ in end $M _l \subset
\R ^3$, $1 \leq k \leq 3$.

The spacetime $({\bf L} ^{3,1}, {\bf g})$ satisfies the dominant
energy condition if, for any timelike vector $W$,
 \bi
 \item[(i)] ${\bf T} _{uv} W ^u W ^v \geq 0$;
 \item[(ii)] ${\bf T} ^{uv} W _u $ is a non-spacelike vector.
 \ei
Restricted on $(M ^3, g, h)$, it implies that
 \beQ
\frac{1}{2}\Big(R +(h ^i _{\;\;i} ) ^2 -h _{ij} h ^{ij}\Big) \geq
\sqrt{\sum _{1\leq i \leq 3} ( \nabla ^j h _{ij} - \nabla _i h ^j
_{\;\;j}) ^2}
 \eeQ
where $R$ is the scalar curvature of $g$.
\begin{thm}
(The Schoen-Yau's positive mass theorem \cite{SY1,SY2, SY3, Wi})
If the spacetime $({\bf L} ^{3,1}, {\bf g})$ satisfies the
dominant energy condition, then, for asymptotically flat initial
data set $(M, g, h)$,
 \beQ
E _l\geq \sqrt{\sum _{1\leq k \leq 3} P _{lk} ^2}
 \eeQ
for each end $M _l$. Equality implies that $M$ has only one end
and ${\bf L} ^{3,1}$ is flat along $M$
\end{thm}

In 1999, the author generalized the positive mass theorem to the
spacetimes including the total angular momentum \cite{Z1}. The
idea is to prove a positive mass theorem for a nonsymmetric
initial data set $(M, g, p)$ where $p$ is a 2-tensor which is not
necessarily symmetric.

\mysection{The positive mass theorem near null infinity} \ls

In spacetimes, null hypersurfaces consist of null geodesics. Along
null hypersurfaces, the induced metrics degenerate and the
geometric properties are difficult to study. However, we can
choose certain spacelike hypersurfaces to approach null infinity
and use them to study null infinity. In Minkowski spacetime, the
spacelike hypersurface
 \beQ
t=\sqrt{1+r ^2}
 \eeQ
has the hyperbolic metric $\breve{g}$ and the nontrivial second
form $\breve{h}$
 \beQ
\breve{g} &=&\frac{ dr ^2 }{1+r ^2}+r ^2 \Big(d \theta ^2 +
\sin ^2 \theta d\psi ^2\Big),\\
 \breve{h} &=&\frac{ dr ^2 }{1+r ^2}+r ^2 \Big(d \theta ^2 +
\sin ^2 \theta d\psi ^2\Big).
 \eeQ
Denote the associated orthonormal frame $\{\breve{e} _i \}$ and
coframe $\{\breve{e} ^i\}$ by $\breve{e} _1 =  \sqrt{1+r
^2}\frac{\pa}{\pa r}$, $\breve{e} _2 =  \frac{1}{r} \frac{\pa}{\pa
\theta}$, $\breve{e} _3 =  \frac{1}{r \sin \theta}\frac{\pa}{\pa
\psi}$, $\breve{e} ^1 = \frac{dr}{\sqrt{1+r ^2} }$, $\breve{e} ^2
= r d \theta $, $\breve{e} ^3 = r \sin \theta d\psi $. Denote
$\breve{\nabla}$ the Levi-Civita connection of $\breve{g}$ and
$\breve{\nabla} _{\breve{e} _i}$ by $\breve{\nabla} _i$. Based on
the above model, we can define asymptotically null initial data
set $(M, g, p)$ ($p$ is not necessarily symmetric) of order $\tau
$ if there is a compact set $K \subset M$ such that $M \setminus
K$ is the disjoint union of a finite number of subsets
$M_1,\cdots, M _l$ - called the ``ends'' of $M$ - each
diffeomorphic to $\R ^3 \setminus B _ r$, where $B _r$ is the
closed ball of radius $r$ with center at the coordinate origin. In
each end, the metric $g$ and the 2-tensor $p$ are $g (\breve{e}_i,
\breve{e}_j)= \breve{g} (\breve{e}_i, \breve{e}_j) +a _{ij}$, $p
(\breve{e}_i,  \breve{e}_j)= \breve{p} (\breve{e}_i, \breve{e}_j)
+b _{ij}$ where $a _{ij}$ and $b _{ij}$ satisfy
 \beQ
 a_{ij}=O\big( \frac{1}{r ^\tau } \big),
 \breve{\nabla} _k  a_{ij}=O\big( \frac{1}{r ^\tau } \big),
 \breve{\nabla} _k \breve{\nabla} _l a_{ij}=
                     O\big( \frac{1}{r ^\tau } \big),
    b _{ij} =O\big( \frac{1}{r ^\tau } \big) ,
    \breve{\nabla} _k  b_{ij}=O\big( \frac{1}{r ^\tau } \big).
 \eeQ

Denote
 \beQ
 \E  &=&
\breve{\nabla} ^ {j} a _{1j} - \breve{\nabla} _1 tr _{\breve{g}}
(a)-\big(a_{11}-g _{11}tr _{\breve{g}} (a) \big), \\
 \P _{k} &=& b_{k1}-g _{k1}tr _{\breve{g}}(b).
 \eeQ
The total energy and the total linear momentum of an
asymptotically null initial data set $(M, g, p)$ on end $M_l$ are
 \beQ
 E_{l\nu} &=& \frac{1}{16\pi}\lim _{r \rw \infty}
\int _{S_{r}} \E  n ^{\nu} r \breve{\omega} _2 \wedge
\breve{\omega} _3,
\label{toe}\\
 P_{l\nu ,k} &=&
\frac{1}{8\pi}\lim _{r \rw \infty} \int _{S_{r}} \P _{k} n ^{\nu}
r \breve{\omega} _2 \wedge \breve{\omega} _3,\label{tolm}
 \eeQ
where $S _{r,l}$ is the sphere of radius $r$ in end $M _l \subset
\R ^3$, $\nu =0,1,2,3$, $k=1,2,3$.
 \begin{thm} \label{pmt-null}
(The positive mass theorem near null infinity \cite{Z3}) Let
 $(M, g, p)$ be a 3-dimensional asymptotically null initial data set of
order $\tau =3$. Let
 \beQ
 \mu &=&
    \frac{1}{2}\big(R+(p _{i} ^{\;i})^2-p _{ij} p ^{ij} \big),\\
 \varpi _j &=& \nabla ^i p _{ji}-\nabla _j p _{i} ^{\;i},\\
 \sigma _j &=& 2\nabla ^i \big( p _{ij}- p _{ji}\big)
 \eeQ
where $R$ is the scalar curvature of $g$. If the initial data set
satisfies the dominant energy condition
 \beQ
 \mu \geq \max \Big\{\sqrt {\sum _{1\leq j \leq 3} \varpi ^2 _j},
 \sqrt {\sum _{1\leq j \leq 3} ( \varpi _j + \sigma _j) ^2}\Big\},
 \eeQ
then,
 \beQ
 E _{l0} - P _{l0,1} \geq \sqrt {\sum _{1 \leq i\leq 3}
 \big(E _{li} - P _{li,1}\big) ^2}.
 \eeQ
If equality holds, then $M$ has only one end and
 \beQ
 R _{ijkl}+p _{ik} p _{jl}-p _{il} p _{jk}=0, \;\;
 \nabla _i p _{jk} - \nabla _j p _{ik}=0, \;\;
 \nabla ^j \big(p _{ij} -p _{ji}\big) =0.
 \eeQ
 \end{thm}
 \begin{rmk}\label{order-tau}
 The proof of Theorem \ref{pmt-null} is also valid for the case $\tau
> \frac{3}{2}$ and the $E _{l\nu} - P _{l\nu,1}$ are finite for $\nu
=0,1,2,3$.
 \end{rmk}
 \begin{rmk}
Theorem \ref{pmt-null} and its application to the positivity of
the Bondi mass (Remark 5.1 in \cite{Z3}) were basically proved in
November 2002 in early version of \cite{Z3}. The final revised
version of \cite{Z3} was sent to Chru\'sciel on June 27, 2003. On
July 23, 2003, the author received the preprint \cite{CJL} from
Chru\'sciel. Translating into our formulation, they proved the
same positive mass theorem as Theorem \ref{pmt-null} for
3-dimensional initial data set $(M, g, K)$ in \cite{CJL}, where
$g$ satisfies integrable condition $(3.7)$ in \cite{CJL}, and $K$
satisfies
 \beQ
K ^{ij} = \frac{1}{r ^3} L ^{ij} +\frac{tr _g K}{3} g ^{ij}, \;\;
tr _g K = C +O\big(\frac{1}{r ^2}\big)
 \eeQ
for certain trace-free tensor $L ^{ij}$ and constant $C$.
 \end{rmk}
\mysection{Bondi's radiating spacetimes}
\ls

Bondi's radiating spacetimes are vacuum spacetimes equipped with
the following metric
 \beq
{\bf g} _{Bondi} &=&\Big(\frac{V}{r} e ^{2\beta}
          +r ^2 e ^{2 \gamma} U ^2 \cosh 2\delta
          +r ^2 e ^{-2 \gamma} W ^2 \cosh 2\delta  \nonumber\\
         & &+2 r ^2 UW \sinh 2 \delta \Big)du ^2
            -2e ^{2\beta} du dr  \nonumber\\
         & &-2r ^2 \Big(e ^{2 \gamma} U \cosh 2\delta
             +W \sinh 2 \delta\Big) du d\theta     \nonumber\\
         & &-2r ^2 \Big(e ^{-2 \gamma} W \cosh 2\delta
            +U \sinh 2\delta \Big)\sin \theta du d\psi    \nonumber\\
         & &+r ^2 \Big(e ^{2 \gamma} \cosh 2\delta d\theta ^2
            +e ^{-2 \gamma}\cosh 2\delta \sin ^2 \theta d \psi ^2  \nonumber\\
         & &+2 \sinh 2\delta \sin \theta d \theta d \psi \Big)  \label{bondi}
 \eeq
where parameters $r>0$, $0 \leq \theta <\pi$, $0 \leq \psi < 2
\pi$ and $\beta, \gamma, \delta, U, V, W$ are smooth functions of
 \beQ
x ^0 =u,\;\;x ^1=r,\;\;x ^2=\theta, \;\;x ^3=\psi.
 \eeQ
The parameter $u$ is called the retarded coordinate. Physically,
$u=constant$ requires to be null hypersurfaces. In Schwarzschild
spacetime, the retarded coordinate $u=t-r-2m \ln \big|r-2m \big|$,
and the metric is written also as
 \beQ
{\bf{g}} _{Sch} =-\Big(1-\frac{2m}{r}\Big)du ^2 -2dudr+r ^2
\Big(d\theta ^2 +\sin ^2 \theta d \psi ^2 \Big).
 \eeQ

The metric (\ref{bondi}) was studied by Bondi, van der Burg,
Metzner and Sachs in the theory of gravitational waves in general
relativity \cite{BBM, Sa, vdB}. They proved that the following
asymptotic behavior holds for $r$ sufficiently large if the
spacetime satisfies the outgoing radiation condition \cite{vdB}
 \beQ
\gamma &=&\frac{c(u, \theta, \psi)}{r} +\frac{C(u,\theta,
\psi)-\frac{1}{6} c ^3 -\frac{3}{2} c d ^2}{r ^3}
+O\big(\frac{1}{r ^4}\big),  \\
\delta &=&\frac{d(u, \theta, \psi)}{r} +\frac{H(u,\theta,
\psi)+\frac{1}{2} c ^2 d -\frac{1}{6} d ^3}{r ^3}
+O\big(\frac{1}{r ^4}\big),  \\
\beta &=&-\frac{c ^2 + d ^2}{4r ^2} +O\big(\frac{1}{r ^4}\big), \\
U &=& -\frac{l(u, \theta, \psi)}{r ^2} +\frac{p(u, \theta,
\psi)}{r ^3}
+O\big(\frac{1}{r ^4}\big),  \\
W &=& -\frac{\bar l(u, \theta, \psi)}{r ^2}+\frac{\bar p(u,
\theta, \psi)}{r ^3}+O\big(\frac{1}{r ^4}\big),  \\
V &=& -r +2 M (u, \theta, \psi) +O\big(\frac{1}{r}\big),
 \eeQ
where
 \beQ
l &=& c _{, 2} +2c \cot \theta +d _{, 3} \csc \theta,\\
\bar l &=& d _{, 2} +2d \cot \theta -c _{,3} \csc \theta,\\
p &=& 2N(u, \theta, \psi) +3(c c _{, 2}+d d _{, 2}) +4(c ^2 +d ^2)\cot \theta\\
  & &-2(c_{,3} d -c d _{,3}) \csc \theta,\\
\bar p &=& 2P(u, \theta, \psi) +2(c _{, 2} d -c d _{, 2}) +3(c
c_{,3} +d d _{,3})\csc \theta.
 \eeQ
Under these conditions, the Bondi's radiating metric ${\bf g}
_{Bondi}$ is
  \beQ
& &-\Big(1-\frac{2M}{r} +O\big(\frac{1}{r ^2}\big)\Big)du ^2
-2\Big(1-\frac{c ^2+d ^2}{4r ^2} +O\big(\frac{1}{r ^4}\big)\Big)du
dr \\
 & &+2 \Big(l+\frac{2c l+2d \bar l}{r} +O\big(\frac{1}{r
^2}\big)\Big)du d\theta +2 \Big(\bar l-\frac{2c\bar l-2d l}{r}
+O\big(\frac{1}{r ^2}\big)\Big)\sin \theta du d\psi \\
 & &+r ^2
\Big(1+\frac{2c}{r} +O\big(\frac{1}{r ^2}\big)\Big)d\theta ^2 +r
^2 \Big(1-\frac{2c}{r}
+O\big(\frac{1}{r ^2}\big)\Big)\sin ^2 \theta d \psi ^2 \\
& &+r ^2 \Big(\frac{4d}{r} +O\big(\frac{1}{r ^2}\big)\Big)\sin
\theta d \theta d \psi.
 \eeQ
Since $\frac{\pa u}{\pa r} \neq 0$ in general, the metric
(\ref{bondi}) is not asymptotically flat at spatial infinity. We
assume
 \begin{description}
 \item[Condition A] Each of the six functions $\beta $, $\gamma$,
$\delta$, $U$, $V$, $W$ and its derivatives up to the second
orders have the same values at $\psi=0$ and $\psi=2\pi$.
 \item[Condition B] For all $u$, $\theta _0 =0$, or $\pi$,
 \beQ
\int _0 ^{2\pi} c(u, \theta _0, \psi) d\psi =0.
 \eeQ
 \end{description}

The Bondi energy-momentum of $u _0$-slice is defined as
 \beQ
m _\nu (u _0) = \frac{1}{4 \pi} \int _{S ^2} M (u _0, \theta,
\psi) n ^{\nu} d S
 \eeQ
for $\nu =0,1,2,3$, where
 \beQ
n ^0 =1,\;\; n ^1 = \sin \theta \cos \psi,\;\; n ^2 = \sin \theta
\sin \psi,\;\; n ^3 = \cos \theta.
 \eeQ
The Bondi energy-momentum is the total energy-momentum measured
after the loss due to the gravitational radiation up to that time.

In 1962, Bondi proved that the $m _{0} (u)$ is a non-increasing
function of $u$ \cite{BBM}, i.e., more and more energy is radiated
away.
 \begin{prop}
(Huang-Yau-Zhang \cite{HYZ}) Let $\big({\bf  L} ^{3,1}, {\bf g}
_{Bondi}\big) $ be a vacuum Bondi's radiating spacetime with
metric ${\bf g} _{Bondi}$ given by (\ref{bondi}). Suppose that
{\bf Condition A} and {\bf Condition B} hold. Then
 \beq
\frac{d}{du} \Big(m _0 - \sqrt{\sum _{1 \leq i \leq 3} m _i ^2 }
\Big) \leq 0. \label{mass-loss}
 \eeq
 \end{prop}
 \pf Denote $|m| =\sqrt{m_1^2+m_2^2+m_3^2 }$. We assume $|m| \neq 0$ otherwise
it reduces to the Bondi mass-loss formula. We have
 \beQ
\frac{d}{du} \Big(m _0 - |m|\Big) &=& -\frac{1}{4\pi}\Big[\int _{S
^2}\Big((c _{,0} )^2 +(d _{,0}
)^2 \Big) d S \\
 & &-\frac{1}{|m|} \sum _{1 \leq i \leq 3}
m_i\,\int _{S ^2}\Big((c _{,0} )^2 +(d _{,0} )^2 \Big) n ^i d S
\Big].
 \eeQ
Using $(n^1)^2+(n^2)^2+(n^3)^2=1$ and H\"older inequality, we
obtain
 \beQ
\sum _{1 \leq i \leq 3}\Big[\int _{S ^2}\Big((c _{,0} )^2 +(d
_{,0} )^2 \Big) n ^i d S \Big] ^2  \leq \Big[\int _{S ^2}\Big((c
_{,0} )^2 +(d _{,0} )^2 \Big) d S \Big] ^2.
 \eeQ
It together with Cauchy-Schwarz inequality implies
 \beQ
\sum _{1 \leq i \leq 3} m_i\,\int _{S ^2}\Big((c _{,0} )^2 +(d
_{,0} )^2 \Big) n ^i d S  \leq  |m| \int _{S ^2}\Big((c _{,0} )^2
+(d _{,0} )^2 \Big) d S.
 \eeQ
Therefore (\ref{mass-loss}) holds. \qed

We study the asymptotically null initial data set $(M, g, h)$
where $M$ is given by
 \beQ
u =\sqrt{1+r ^2} -r + \frac{\big(c ^2 +d ^2\big) _{u=0}}{12r ^3}
+\frac{a _{3} (\theta, \psi)}{r ^4}+a _{4}
 \eeQ
where $a _{4} (r,\theta, \psi)$ is a smooth function which
satisfies: In the Euclidean coordinate systems $\{z ^i\}$, $|z|
=r$, $a_4 =o\big(\frac{1}{r ^4}\big)$, $\pa a_4 =o\big(\frac{1}{r
^5}\big)$, $\pa \pa a_4 =o\big(\frac{1}{r ^6}\big)$ as $r
\rightarrow \infty$. We compute the induced metric $g$ and the
second fundamental form $h$. Define $a \approx b$ if and only if
$a =b +o\big(\frac{1}{r^3}\big)$.
 \beQ
g(\breve{e} _1, \breve{e} _1) & \approx &1+\frac{16a_3+M -c
c_{,0}- d d
_{,0}}{2r^3} \Big | _{u=0},\\
 g(\breve{e} _1, \breve{e} _2)
&\approx &-\frac{l}{2r^2}
 +\frac{12 N -3 l _{,0}
 +4 (c c _{,2}+ d d _{,2})}{12r ^3},\\
g(\breve{e} _1, \breve{e} _3)&\approx&-\frac{\bar l}{2r ^2}
+\frac{12 P-3\bar{l} _{,0}+4\csc \theta (c c _{,3}+  d d _{,3})}{12r ^3},\\
g(\breve{e} _2, \breve{e} _2)&\approx&1+\frac{2c}{r}+\frac{2(c^2+ d^2) +c _{,0}}{r ^2}\\
  & & +\frac{c ^3 +c d ^2 +2C+2(c c _{,0} + d d_{,0})
         +\frac{c _{,00}}{4}}{r ^3},\\
g(\breve{e} _2, \breve{e} _3) &\approx
&\frac{2d}{r}+\frac{d_{,0}}{r^2}
   +\frac{c ^2 d +d ^3 +2H+\frac{d _{,00}}{4}}{r ^3}\\
g(\breve{e} _3, \breve{e} _3) &\approx &1-\frac{2c}{r}+\frac{2(c^2
+d^2) -c _{,0}}{r ^2}\\
     & &+\frac{-c ^3 -c d ^2-2C +2(c c _{,0} +d d _{,0}) -\frac{c
_{,00}}{4}}{r ^3},\\
 h(\breve{e} _1,
\breve{e} _1)
   &\approx&1 + \frac{c ^2+d^2}{r ^2}
   +\frac{16a_3-M}{r ^3},\\
h(\breve{e} _1, \breve{e} _2)
   &\approx & \frac{l}{2r ^2}
      +\frac{1}{2r^3}\big[ \frac{l _{,0}}{2}-2 (c ^2+ d ^2)  \cot\theta -4 N\\
   & & (-c d _{,3} +c _{,3}d )\csc \theta -\frac{13}{3}(c c _{,2} +d d _{,2})\big],\\
 h(\breve{e} _1, \breve{e} _3)
  &\approx & \frac{\bar l}{2r^2 }
      +\frac{1}{2r ^3 }\big[\frac{\bar l _{,0}}{2} +c d _{,2}-c _{,2} d
      -4P\\
   & & -\frac{13}{3} (c c _{,3} +d d _{,3}) \csc \theta\big],\\
h(\breve{e} _2, \breve{e} _2)&\approx &1 +\frac{c}{r}
     +\frac{c _{,0}}{r ^2}
   +\frac{1}{4r ^3}\big[ 3M-16 a_3 -4C -2 l _{,2} \\
   & &- 2c(c ^2 + d ^2) +5(c c _{,0} +d d _{,0})+\frac{3}{2}c
   _{,00}\big],\\
h(\breve{e} _2, \breve{e} _3)
       &\approx&\frac{d}{r}+\frac{d _{,0}}{r^2}
        +\frac{1}{4r ^3 }\big[-2d(c ^2 +d ^2) +2d \cot ^2 \theta
        \\
& &+2d \csc ^2 \theta -4c _{,3} \cot \theta \csc \theta -d
_{,33}\csc ^2 \theta
 \eeQ
 \beQ
 & & -d _{,2}\cot \theta - d _{,22}
          - 4H +\frac{3}{2} d_{,00} \big],\\
h(\breve{e} _3, \breve{e} _3)&\approx &1 -\frac{c}{r}
     -\frac{c _{,0}}{r ^2}
   +\frac{1}{4r ^3}\big[ 3M-16a_3 +4C \\
   & &+2c(c^2 + d ^2) +5(cc _{,0}
   +dd _{,0})-\frac{3}{2}c _{,00}\\
   & &-2 l \cot \theta -2 \bar l _{,3} \csc \theta \big].
 \eeQ
Here all functions in the right hand sides take value at $u=0$ and
all derivatives with respect to $x ^2$ and $x ^3$ are taken after
substituting $u=0$.

Replacing $u$ by $u-u _0$, applying Theorem \ref{pmt-null}, Remark
\ref{order-tau} to $(M, g, h)$ and using (\ref{mass-loss}), we can
prove the following theorem concerning positivity of the Bondi
mass.
 \begin{thm}
(Huang-Yau-Zhang \cite{HYZ}) Let $\big({\bf  L} ^{3,1}, {\bf g}
_{Bondi}\big) $ be a vacuum Bondi's radiating spacetime with
metric ${\bf g} _{Bondi}$ given by (\ref{bondi}). Suppose that
{\bf Condition A} and {\bf Condition B} hold. If there is $u _0$
such that $c \big| _{u=u_0} = d \big| _{u=u_0} = 0$, then
 \beQ
m _{0} (u) \geq \sqrt{\sum _{1 \leq i \leq 3} m ^2 _{i} (u)}
 \eeQ
for all $u \leq u_0$. Moreover, if the equality holds for all $u
\leq u_0$, ${\bf L} ^{3,1}$ is flat in the region $u \leq u_0$.
 \end{thm}

It turns out the Bondi mass cannot become negative if
 \beQ
\lim _{u \rightarrow \infty} c(u, \theta, \psi)=\lim _{u
\rightarrow \infty} d(u, \theta, \psi)=0.
 \eeQ
In \cite{HYZ}, we are also studying the positivity of the Bondi
mass by using Schoen-Yau's argument in \cite{SY4}.

\end{document}